\documentclass[12pt]{amsart}

\usepackage[dvipsnames]{xcolor}
\usepackage[shortlabels]{enumitem}
\usepackage[pagebackref,colorlinks,linkcolor=Maroon,citecolor=MidnightBlue,urlcolor=NavyBlue,hypertexnames=true]{hyperref}
\usepackage{xurl}

\newtheorem{theorem}{Theorem}[section]
\newtheorem{proposition}[theorem]{Proposition}

\newtheorem{lemma}[theorem]{Lemma}
\theoremstyle{definition}

\theoremstyle{remark}
\newtheorem{remark}[theorem]{Remark}

\newcommand{\R}{\mathbb R}
 \newcommand{\C}{\mathbb C}
\newcommand{\A}{\mathcal A}
 \newcommand{\B}{\mathcal B}
 \newcommand{\K}{\mathcal K}
\newcommand{\Hh}{\mathcal H}
\newcommand{\Hl}{\mathcal L}

\title[Equivalence of Linear $\&$ C* extremality for Real UCP maps]{On the equivalence of Linear and $C^*$ extremality for unital Completely Positive Maps on commutative Real $C^*$-algebras}
\author[S. Balasubramanian]{Sriram Balasubramanian}
\address{Department of Mathematics\\
IIT Madras, Chennai - 600036, India.}
\email{bsriram80@yahoo.co.in, bsriram@iitm.ac.in}

\author[T. Karthikeyan]{Karthikeyan Thangaraj}
\address{Department of Mathematics\\
IIT Madras, Chennai - 600036, India.}
\email{ma24d203@smail.iitm.ac.in, kk.mathsbiology123@gmail.com}



\begin{document}
\maketitle

\begin{abstract}
Let $\A$ be a unital commutative real $C^*$-algebra and 
$\Hh$ be a finite dimensional real Hilbert space. 
In this article, an easily verifiable necessary and sufficient condition 
for a linear extreme unital completely positive (UCP) map 
$\Phi:\A \rightarrow B(\Hh)$ 
to be C*-extreme, is obtained.  As an application, the 
$C^*$-algebras $\A$ on which linear extreme UCP maps coincide 
with the $C^*$-extreme ones, are classified. 
\end{abstract}


\section{Introduction}
Convex analysis plays a vital role in understanding 
various important classes of positive maps between $C^*$-algebras, 
such as completely positive maps, PPT maps, and entanglement breaking maps.
 Our main goal here is to better understand the convex extremal structure 
 of unital completely positive maps in the real setting. \\

The study of completely positive maps on complex $C^*$-algebras 
was initiated by Stinespring and Arveson in the seminal papers \cite{S}, \cite{A}. 
Among various significant results in \cite{A}, Arveson obtained 
 an abstract characterization of the linear
extreme points of the convex set of completely positive maps from  
a complex $C^*$-algebra to $B(\K)$ for some Hilbert space $\K$,
 in terms of minimal Stinespring dilations.  A noncommutative notion 
of convexity that has received considerable attention is $C^*$-convexity, which 
was introduced and studied by Loebl, Moore and Paulsen in \cite{LP} and
\cite{HMP}.  There is extensive literature available on linear and 
$C^*$-extreme points of different classes of completely positive maps between complex
 $C^*$-algebras. See for example, \cite{A}, \cite{BBK}, \cite{BK}, \cite{BDMS}, 
\cite{BH},  \cite{F}, \cite{FM1}, \cite{FM2}, \cite{FZ}, \cite{G1}, \cite{G2}, 
\cite{Z} and the references therein. \\

 In this paper, we work with the 
convex set $\operatorname{UCP}(\A, B(\Hh))$ of  all unital completely positive maps from 
$\A$ to $B(\Hh)$, where $\A$ is a unital commutative real $C^*$-algebra and $\Hh$ 
is a finite-dimensional real Hilbert space. This work is motivated by 
 \cite{ANSS}, in which the authors obtained a complete description 
of the $C^*$-extreme points of this set. Please see \cite[Theorem 4.10]{ANSS}. 
In particular, they also observed that every $C^*$-extreme point of 
$\operatorname{UCP}(\A, B(\Hh))$ is also a linear extreme point, i.e. an extreme point in 
the usual convexity sense. \\

There are two main contributions in this article. Firstly, 
we obtain a novel characterization of the linear extreme points in 
$\operatorname{UCP}(\A, B(\Hh))$ that are also $C^*$-extreme. Secondly, as an 
application, we characterize all the real unital commutative $C^*$-algebras $\A$ for 
which the linear extreme points of $\operatorname{UCP}(\A, B(\Hh))$ coincide 
with the $C^*$-extreme points. The required definitions are given below. 
Several basic facts about CP maps that we have used here have been drawn from \cite{P}. \\

First, 
we characterize the linear extreme points in this set  
that are also $C^*$-extreme. Second, as an application, we characterize all such 
$C^*$-algebras $\A$ on which linear extreme points coincide with the
 $C^*$-extreme points of this set. The required definitions are given below. 
 Several basic facts about CP maps that we have  
used here have been drawn from \cite{P}. \\


We retain most of the notations 
from \cite{ANSS} and begin by first recalling the well-known result 
(see \cite[Proposition 5.1.4]{L}) that $\A$ is (isometrically) *-isomorphic to the real commutative $C^*$-algebra 
\[
C(\Omega,-)
=\{f\in C(\Omega,\C):f(\bar w)=\overline{f(w)}\},
\]
where $\Omega$ is a compact Hausdorff topological space 
and $-$ is a homeomorphism of order at most two on $\Omega$ 
 that maps $w \mapsto\bar w$.  Note that the involution
 $*$ on $C(\Omega,-)$ is given by $f^*(w) = \overline{f(w)}$. 
Under this identification, 
throughout the paper we will simply assume $\A = C(\Omega,-)$. 
The standard complexification of $\A$, denoted by $\A_c$, will 
also be invoked wherever necessary.  \\

A point $\xi \in \Omega$ will be called  
{\bf symmetric} if $\xi=\bar\xi$, and {\bf non-symmetric} if 
$\xi \neq \bar\xi$. 
For every $w\in\Omega$, define the one-dimensional UCP map $\rho_w$ by 
\begin{equation}
\label{eq:rho-w}
\rho_w(f)=Re f(w),\qquad f\in C(\Omega,-),
\end{equation}
and the two-dimensional UCP map $\Pi_w$ by 
\begin{equation}
\label{eq:Pi-w}
\Pi_w(f)=
\begin{pmatrix}
Re f(w)& Im f(w)\\
-Im f(w)& Re f(w)
\end{pmatrix}.
\end{equation}
If $w=\bar w$, then $f(w)\in\R$ for every $f\in C(\Omega,-)$, so 
$\rho_w(f)=f(w)$ is the one-dimensional real irreducible representation
associated with the symmetric point $w$. If $w\ne\bar w$, then $\Pi_w$ is 
the two-dimensional irreducible real representation associated with the 
non-symmetric point $w$, while $\rho_w$ is its 
one-dimensional compression. Note that $\rho_w$, in this case, is not 
multiplicative. \\

Let $\mathcal B$ and $\mathcal B'$ denote unital real $C^*$-algebras.
A (real) linear map
$
\Phi:\mathcal B\rightarrow \mathcal B'
$
is said to be {\bf unital} if
$
\Phi(1_{\mathcal B})= 1_{\mathcal B'}  
$
and {\bf positive} if
$
\Phi(a)\ge 0,
$
in $\mathcal B'$, 
whenever $a \ge 0$ in $\mathcal B$.
For each $k\in\mathbb N$, the $k$-th amplification of $\Phi$ is the map
$
\Phi_k:M_k(\mathcal B)\rightarrow M_k({\mathcal B'})
$
defined by
$
\Phi_k([a_{ij}])=[\Phi(a_{ij})].
$
The map $\Phi$ is said to be {\bf completely positive} if
$\Phi_k$ is positive for every $k\in\mathbb N$. A unital completely
positive map will be called a {\bf UCP} map. Throughout the paper, the 
notation 
$
\operatorname{UCP}(\mathcal B,{\mathcal B'})
$
will be used to denote the (convex) set of all UCP maps from 
$\mathcal B$ into ${\mathcal B'}$. \\

Next we define the two necessary notions of extreme points. A map  
$
\Phi\in \operatorname{UCP}(\mathcal B,{\mathcal B'}) 
$
is said to be a {\bf linear extreme point} of 
$\operatorname{UCP}(\mathcal B,{\mathcal B'})$ if 
$
\Phi=\sum_{i=1}^r t_i\Phi_i,
$
for some 
$
\Phi_i\in \operatorname{UCP}(\mathcal B,{\mathcal B'})$,
$t_i\in(0,1)$ satisfying 
$\sum_{i=1}^r t_i=1,
$
then 
$
\Phi_i=\Phi
$
for each $i$. \\

Given 
$
\Phi_1,\ldots,\Phi_r
\in
\operatorname{UCP}(\mathcal B,{\mathcal B'}), 
$
the linear map 
$
\Phi\in \operatorname{UCP}(\mathcal B,{\mathcal B'})
$
is said to be a {\bf $C^*$-convex combination} of the maps
$\Phi_1,\ldots,\Phi_r$ if there exist 
$
x_1,\ldots,x_r\in {\mathcal B'}
$
such that
$
\sum_{i=1}^r x_i^*x_i=1_{\mathcal B'}
$
and
$
\Phi(a)
=
\sum_{i=1}^r
x_i^*\Phi_i(a)x_i,
$
for all $a\in\mathcal B.$ 
The $C^*$-convex combination is called {\bf proper} if each 
$x_i$ is invertible. It is called {\bf trivial} if 
for each $i$, there exists
a unitary element $u_i\in {\mathcal B'}$ such that 
$
\Phi_i(a)=u_i^*\Phi(a)u_i,
$
for all $a \in\mathcal B.$ A linear map 
$
\Phi\in \operatorname{UCP}(\mathcal B,{\mathcal B'})
$
is said to be a {\bf $C^*$-extreme point} of 
$\operatorname{UCP}(\mathcal B, {\mathcal B'})$, if every 
expression of $\Phi$ as a proper $C^*$-convex combination is trivial. 
Indeed, linear extremality coincides with $C^*$-extremality, when 
${\mathcal B'}$ is one dimensional. \\

Recall the UCP maps $\rho_w$ and $\Pi_w$ from equations 
\eqref{eq:rho-w} and \eqref{eq:Pi-w}.    Henceforth, if $r \ge 1$, 
the notation $I_r$ will denote the $r \times r$ identity matrix, the 
operator $A^{(r)}$ will denote the direct sum of $r$ copies of the operator $A$. 
And, in a direct sum with $k$ summands, whenever $k=0$, it is to be  
understood that the direct sum is absent. \\

The following is a structure theorem  for $C^*$-extreme UCP maps from
 $\A \to B(\Hh)$. 

\begin{theorem}\cite[Theorem~4.10]{ANSS}
\label{thm:ANSS-main}
Let $\A=C(\Omega,-)$. A linear map 
$\Phi \in \operatorname{UCP}(\mathcal A,B(\Hh))$ 
is $C^*$-extreme 
if and only if, 
\begin{equation}
\label{eq:structure}
\Phi \sim 
\left( 
\bigoplus_{i=1}^{p}
\rho_{\xi_i}^{(m_i)}
\right)
\oplus 
\left(
\bigoplus_{j=1}^{q}
  \rho_{\gamma_j}
\right)
\oplus
\left(
\bigoplus_{k=1}^{r}
\Pi_{\eta_k}^{(s_k)} 
\right),
\end{equation}
where '$\sim$' denotes unitary equivalence, $p,q, r \in 
\{0,1,2,...\}$ but not simultaneously zero, 
$\xi_i$'s are distinct symmetric points,  
$\gamma_j, \eta_k$ are non-symmetric points 
satisfying $\{\gamma_s,\bar\gamma_s\} 
\neq \{\gamma_t,\bar\gamma_t\}$ and $\{\eta_s,\bar\eta_s\} 
\neq \{\eta_t,\bar\eta_t\}$ when     
$s \neq t$, $m_i, s_k \in \{1,2,\dots\}$ and $\sum_i m_i + 
q + \sum_k 2s_k = dim(\Hh)$. 
\end{theorem}


\begin{remark}
    \label{rem:oneway}
That $
\Phi\text{ is } C^*\text{-extreme}
\implies
\Phi\text{ is linear extreme},$
follows from \cite[Proposition~3.8]{ANSS}.
\end{remark}

Recall that $\Hh$ is a finite dimensional real Hilbert space. 
Our main results are the following.
\begin{theorem}
\label{thm:range}
Let $\A=C(\Omega,-)$ and $\Phi \in \operatorname{UCP}(\mathcal A,B(\mathcal H))$. 
The following statements are equivalent.
\begin{enumerate}[\rm(i)]
\item $\Phi$ is $C^*$-extreme.
\item $\Phi$ is linear extreme with commutative range. 
\end{enumerate}
\end{theorem}

A proof of Theorem \ref{thm:range} is given in Section
\ref{sec:main1}.  As an application, we obtain: 

\begin{theorem}
\label{thm:dom}
Let $\A=C(\Omega,-)$ and $dim(\Hh) \ge 2$. 
 The following statements are equivalent. 
\begin{enumerate}[\rm(i)]
\item Every linear extreme 
$\Phi \in \operatorname{UCP}(\mathcal A,B(\mathcal H))$ is $C^*$-extreme.
\item $\A\cong\R,\, \R^2\text{ or }\C.$
\item $\dim(\A) \le 2.$
\end{enumerate}
\end{theorem}

A proof of Theorem \ref{thm:dom} can be found in Section \ref{sec:main2}. 

\section{Preliminaries} 
\label{sec:prelims}
In this section, we include some well-known results and also record some important
 observations. While some of the results stated 
here hold on a much higher 
level of generality, we only state versions directly applicable to our context, 
starting with a very useful characterization of real UCP maps. 

\begin{theorem} \cite[Theorem 1]{S}
\label{thm:stinespring}
Let $\mathcal B$ be a unital real $C^*$-algebra and 
$\Hh$ be a finite-dimensional real Hilbert space. The following 
statements are equivalent.
\begin{enumerate}[\rm(i)]
\item $\Psi \in \operatorname{UCP}(\mathcal B,B(\mathcal H))$. 
\item There exists a real Hilbert space $\mathcal K$, 
an isometry $V: \Hh \rightarrow \mathcal K$ 
and a $*$-homomorphism $\pi:\mathcal B \rightarrow B(\mathcal K)$ such that $\Psi(a) = V^* \pi(a) V$. 
\end{enumerate}
\end{theorem} 

Note that \cite[Theorem 1]{S} in its original form 
is the complex version of the above stated theorem.  But the same proof works fine in the real setting as well, as was also pointed out in \cite{ANSS}. 
The triple $(\mathcal K, \pi, V)$ in 
Theorem \ref{thm:stinespring} is 
called a Stinespring representation of $\Psi$. If $\mathcal K = \overline{span}\{\pi(a)Vh\,:\,a \in \mathcal B, \, h \in \Hh\}$, then 
the triple $(\mathcal K, \pi, V)$ is called a {\bf minimal Stinespring representation of $\Psi$}. \\





We next recall the Arveson criterion to determine 
linear extremality of the UCP map $\Psi$ in terms of its minimal Stinespring 
representation.

\begin{theorem}\cite[Theorem 1.4.6]{A} 
\label{thm:arveson}
Let $\mathcal B$ be a unital real $C^*$-algebra and 
$\Hh$ be a finite-dimensional real Hilbert space. Suppose that 
$\Psi \in \operatorname{UCP}(\mathcal B,B(\mathcal H))$ and
 $(\K,\pi,V)$ be its minimal Stinespring representation. 
The following statements are equivalent.
\begin{enumerate}[\rm(i)]
\item $\Psi$ is linear extreme. 
\item $T=T^*\in\pi(\A)'$ and $V^*TV=0$ $\implies$ $T=0$.
\end{enumerate}
\end{theorem}

Note that \cite[Theorem 1.4.6]{A} in its original form 
is the complex version of the above stated theorem. 
But the same proof works fine in the real setting as well, 
as was also pointed out in \cite{ANSS}. \\

Next we prove two important lemmas that will be used in the sequel.

\begin{lemma}
\label{lem:dirsum} 
Let $\mathcal B$ be a real unital $C^*$-algebra, $\mathcal E$ and $\Hl$ be real Hilbert spaces, 
$\Psi \in \operatorname{UCP}(\mathcal B,B(\mathcal E))$ and 
$\rho:\mathcal B \to B(\Hl)$ is a unital $*$-homomorphism. 
The following statements are equivalent.
\begin{enumerate}[\rm(i)]
\item $\Psi \in \operatorname{UCP}(\mathcal B,B(\mathcal E))$ is linear extreme. 
\item $\Psi \oplus\rho:\A\rightarrow B(\mathcal E \oplus \Hl)$ 
is linear extreme. 
\end{enumerate}
\end{lemma}

\begin{proof} 
$(ii) \implies (i):$ If $\Psi = \frac 1 2\{\Psi_1 + \Psi_2\}$ for some UCP maps, then $\Psi \oplus \rho = 
\frac 12 \{(\Psi_1 \oplus \rho) + (\Psi_2 \oplus \rho)\}$. It follows from the hypothesis that 
$\Psi_1 = \Psi_2$. \\

$(i) \implies (ii):$ Let $(\K,\pi,V)$ be a minimal Stinespring representation of $\Psi$. Recall that 
$
\Psi(a)=V^*\pi(a)V,
$
for all $a \in \mathcal B$. 
Since $\rho$ is already a $*$-representation, a Stinespring
representation of $\Psi\oplus\rho$ is given by 
$
(\K\oplus {\Hl},\,\pi\oplus\rho,\,V\oplus I_{\Hl}).
$
The minimality of this Stinespring representation is a direct 
consequence of the minimality of the Stinespring representation 
$(\K, \pi, V)$ of $\Psi$. \\

To establish the linear extremality of $\Psi \oplus \rho$, we use Theorem 
\ref{thm:arveson}.  Let
$
T= \begin{pmatrix}
A&B\\
B^*&C
\end{pmatrix} 
 \in ((\pi\oplus\rho)(\mathcal B))'. 
$
Because $T$ commutes with $(\pi\oplus\rho)(\mathcal B)$, it must 
be the case that $A \in\pi(\mathcal B)'$, $C\in\rho(\mathcal B)'$, 
and
$
\pi(a)B=B\rho(a), 
$
for all $a\in \mathcal B$. 
Suppose that 
\[
0 = (V\oplus I_{\Hl})^*T(V\oplus I_{\Hl}) = \begin{pmatrix}
V^*AV&V^*B\\
B^*V&C
\end{pmatrix}
\]
Consequently, 
\[
C=0,\quad B^*V=0,\quad V^*B=0,\quad V^*AV=0.
\]
We now show that $B=0$.  Indeed, for every $a\in \mathcal B$ and $h\in \mathcal E$,
\[
B^*\pi(a)Vh=\rho(a)B^*Vh=0.
\]
By minimality of the Stinespring representation $(\K,\pi,V)$, i.e. the denseness of
$ \operatorname{span}\{
\pi(a)Vh \,:\, h \in \mathcal E\}$ in $\mathcal K$, it follows that $B^*=0$.  
Finally, recall that $A=A^* \in\pi(\mathcal B)'$ with $V^*AV=0.$ 
The conclusion $A = 0$ follows from the linear extremality of 
$\Psi$ and an application of Theorem \ref{thm:arveson}.  Thus $T=0$.  
An application of Theorem \ref{thm:arveson} 
completes the proof. 
\end{proof}

\begin{lemma}
\label{lem:comp}
Let $\rho:\mathcal C\rightarrow\B$ be a surjective unital $*$-homomorphism 
of real unital $C^*$-algebras, let $\Psi \in \operatorname{UCP}(\mathcal B,M_n(\R))$ and let 
$\Phi:=\Psi \circ \rho \in \operatorname{UCP}(\mathcal C, M_n(\R))$. The following statements are equivalent.
\begin{enumerate}[\rm(i)]
\item $\Psi \in \operatorname{UCP}(\mathcal B,M_n(\R))$ is linear extreme.
\item $\Phi \in \operatorname{UCP}(\mathcal C, M_n(\R))$ is linear extreme.  
\end{enumerate}
\end{lemma}

\begin{proof}
$(i) \implies (ii):$ Let $(\mathcal K, \pi,V)$ be a minimal Stinespring 
representation for the UCP map $\Psi$. Since $\Phi = \Psi \circ \rho$, for 
each $x \in \mathcal C$, one gets $\Phi(x) = \Psi (\rho(x)) = V^* (\pi \circ \rho) (x) V$. 
Thus $(\mathcal K, \pi \circ \rho, V)$ is a Stinespring representation of $\Phi$. 
Since $\rho$ is surjective, the range of $\pi$ coincides with the range of $\pi \circ \rho$. 
This implies that $(\mathcal K, \pi \circ \rho, V)$ is in fact, a minimal Stinespring dilation 
for $\Phi$. Note that since the ranges of $\pi$ and $\pi \circ \rho$ are identical, 
so are their commutants.  If $T = T^* \in (\pi \circ \rho)(\mathcal C)' = \pi(\mathcal B)'$ 
satisfies $V^*TV = 0$, then it follows from Theorem \ref{thm:arveson} and the linear extremality 
of $\Psi$, that $T = 0$. Another application of Theorem \ref{thm:arveson}
implies that $\Phi$ is a linear extreme UCP map.  \\

$(ii) \implies (i):$ Suppose that
 $\Psi = \frac 1 2(\Psi_1 + \Psi_2)$ for some UCP maps $\Psi_1 \neq \Psi_2.$  
Let $\Phi_i = \Psi_i \circ \rho$. Observe that $\Phi = 
\frac 1 2 (\Phi_1 + \Phi_2)$. Also, $\Phi_1 \neq \Phi_2$, since $\Psi_1 \neq 
\Psi_2$ and $\rho$ is surjective. This contradicts the hypothesis that the UCP map 
$\Phi$ is linear extreme. 
\end{proof}

\section{Range Classification}
\label{sec:main1}
This section contains a proof of the first main result of the paper, i.e. Theorem 
\ref{thm:range}. We first prove a lemma which is a key tool.   \\

Recall the representation $\Pi_w$ from equation
 \eqref{eq:Pi-w}. It is customary to view $\Pi_w$ as a map  
 from $\A$ to $\C$, where $\C$ is the real 
$C^*$-algebra identified via the injective $*$-homomorphism 
$\Pi:\C \rightarrow M_2(\R)$ given by 
\begin{equation}
\label{eq:Pi}
\Pi(a+ib) = \begin{pmatrix}  
a & b \\ -b & a \end{pmatrix} \in M_2(\R). 
\end{equation}

\begin{lemma}
\label{lem:ext-surj}
Let $\A=C(\Omega,-)$. If 
$
\psi:\A\rightarrow\C
$
is a linear extreme surjective UCP map, then there is a non-symmetric point 
$\eta\in\Omega$ such that, 
\[
\psi(f)=f(\eta),\qquad f\in\A.
\]
Equivalently, 
$\Pi \circ \psi$ is the irreducible representation $\Pi_\eta$. 
\end{lemma}

\begin{proof}
The proof is adapted from \cite[Proposition 4.5]{ANSS}. 
Let $\A_c$ denote the complexification of $\A$.  
It is well-known (see \cite{L}) and also pointed out in \cite{ANSS} that $\A_c = C(\Omega)$, 
the complexification $\C_c$ of the real $C^*$-algebra $\C$ is isomorphic to 
$\C\oplus\C$ via $\Gamma$ defined by, 
\[
\Gamma(z+iw)=\bigl(z+iw,\overline z+i\overline w\bigr),
\qquad z,w\in\C, 
\]
and the complexification of the UCP map $\psi$ is the
complex linear UCP map 
$
\psi_c:\A_c\to \C\oplus\C
$
given by 
\[
\psi_c(f +ig)
 =\bigl(\psi(f)+i\psi(g),\,
        \overline{\psi(f)}+i\,\overline{\psi(g)}\bigr),
\]
for all $f,g \in\A$. 
Let $\tilde{\psi} = V^*\psi_c V$, where $V^* = \begin{pmatrix}  I & 0 \end{pmatrix}:\C^2 \rightarrow \C.$ i.e., 
\[
\tilde{\psi}(f+ig)=\psi(f)+i\psi(g).
\]
Observe that $\tilde{\psi}:\A_c \rightarrow \C$ is a complex linear UCP map. 
Suppose $\tilde{\psi} = \frac12 (\gamma + \delta)$ for some complex linear UCP maps 
$\gamma, \delta:\A_c \rightarrow \C$. Note that
$\gamma$, $\delta$ are real linear UCP maps when restricted to $\A$. Moreover, 
restricting $\tilde{\psi}$ to $\A$ yields, 
\[
\tilde{\psi}(f) = \psi(f) = \frac1 2 (\gamma(f) + \delta(f)),
\]
for each $f \in \A$. By linear extremality of $\psi$, it follows that 
$\gamma(f) = \delta(f) = \psi(f)$ for all $f \in \A$.  By complex linearity of 
$\gamma$, it follows that 
\[
\gamma(f + ig) = \gamma(f) + i \gamma(g) = \psi(f) + i \psi(g) = \tilde{\psi}(f+ig),
 \]
for all $f, g \in \A$. Thus $\tilde{\psi} = \gamma = \delta$ and  
 $\tilde{\psi}$ is a linear extreme complex UCP map.  %
It is well known that linear extreme complex UCP functionals 
on the commutative complex $C^*$-algebra $C(\Omega)$ 
are precisely the point evaluation maps. 
Thus there exists $\eta\in\Omega$ such that 
$
\tilde{\psi}(h)=h(\eta), 
$ for all $h = f + i g \in \A_c$.
Restricting to $\A$ implies 
$\psi(f)=f(\eta),
$ for all $ f\in\A.$ 
If $\eta=\bar\eta$, then for any $f \in \A$,  
$f(\eta) = f(\bar\eta) = \overline{f(\eta)}$. In particular,  
 $f(\eta)\in\R$. This contradicts the surjectivity of $\psi$. 
 Hence $\eta\ne\bar\eta$ and $\Pi \circ \psi = \Pi_\eta$. This completes 
the proof.    
\end{proof}

\subsection{Proof of Theorem \ref{thm:range}}
We are now in a position to prove the
first main result of the paper, i.e. Theorem \ref{thm:range}.
\begin{proof}
$(i) \implies (ii)$: By Proposition \cite[Proposion 3.8]{ANSS}, $\Phi$ 
is linear extreme. Recall the structure of $\Phi$ from equation 
\eqref{eq:structure} in 
Theorem \ref{thm:ANSS-main}. 
Since the UCP maps $\rho_{\xi_i}$, $\rho_{\eta_j}$ and $\Pi_{\eta_j}$ that make up the blocks of 
$\Phi$ have commutative range, so does $\Phi$. \\

To prove $(ii) \implies (i)$, let $\B=C^*(\operatorname{Range}(\Phi))\subset M_n(\R)$. 
Observe that $\B$ is a finite-dimensional unital commutative real $C^*$-algebra. 
It is well known that $\B \cong \R^\ell\oplus\C^m,$ for some $\ell, m \ge 0$. 
Recall $\Pi$ from equation \eqref{eq:Pi}.   
Identifying $\B$ as a $C^*$-subalgebra of  
$M_n(\R)$, by \cite[Theorem 2.5.15]{HJ}, there exist UCP maps    
$
\phi_j:\A\to\R,
\,
\psi_k:\A\to\C
$
such that 
\begin{equation}
\label{eq:Phi-coords}
\Phi(f)=
\bigl(
\bigoplus_{j=1}^\ell
(\phi_j(f))^{(m_j)}
\bigr)
\oplus
\bigl(
\bigoplus_{k=1}^m
\left(\Pi(\psi_k(f)) \right)^{(s_k)}\bigr). 
\end{equation}
Note that if $\ell$ or $m = 0$, then the corresponding direct sum 
in equation \eqref{eq:Phi-coords} is absent. Also, 
$\sum_j m_j + \sum_k 2s_k = n$. 
Since $\Phi$ is linear extreme and $\Pi$ is an injective $*$-homomorphism, 
it follows that the maps $\phi_j$ and $\psi_k$ in equation \eqref{eq:Phi-coords}  
are linear extreme UCP maps as well. Otherwise a non-trivial convex     
decomposition of one of these maps, with all the other maps in equation     
\eqref{eq:Phi-coords} left unchanged, gives a non-trivial convex    
decomposition of $\Phi$ contradicting its linear extremality.   
It follows from \cite[Corollary 4.6(i)]{ANSS} that
$\phi_j=\rho_{\xi_j}$ for some $\xi_j\in\Omega$. Note that $\xi_j$ does not 
necessarily have to be a symmetric point. And if $\xi_j$ is 
non-symmetric, i.e. $\xi_j \neq \bar \xi_j$, 
then $m_j \le 1$. Otherwise, 
\[
\phi_j \oplus \phi_j = \rho_{\xi_j} \oplus \rho_{\xi_j} = \frac 1 2 \left(\Pi_{\xi_j}+ \Pi_{\bar \xi_j} \right),
\]
which implies that $\phi_j \oplus \phi_j$ is a linear extreme 
UCP map. This in turn contradicts the linear extremality  
of $\Phi$. \\


Fix an index $k$ from equation \eqref{eq:Phi-coords}. Observe that $\Pi \circ \psi_k$ 
is surjective. This is because, the corresponding $k$-th coordinate projection map from 
$\B$ to $C^*(\psi_k(\A))$ is surjective.  And since $\Pi$ is injective and 
$\psi$ is unital, it must be the case that $\psi_k$ is surjective. Note that as observed 
earlier, $\psi_k$ is also linear extreme.  Hence, by Lemma \ref{lem:ext-surj}, $\Pi \circ \psi_k = \Pi_{\eta_k}$ for 
some non-symmetric point $\eta_k \in \Omega$, i.e., $\eta_k\ne\bar\eta_k$. \\

Since $\Pi_{\eta_k} \sim \Pi_{\bar \eta_k}$, after a rearrangement of the blocks, (and using 
the fact that corresponding to each $\eta_k$, there is at most one scalar block counting multiplicities)
 it follows that upto unitary equivalence, the UCP map $\Phi$ has exactly the same block form as in 
\cite[Theorem~4.10]{ANSS}. Thus $\Phi$ is $C^*$-extreme, and the proof is 
complete.  


\end{proof}



%

\begin{remark}
\label{rem:range-not-an-alg}
Recall that $\A = C(\Omega, -)$. 
 \vspace{0.2 cm}
 
(i) Let $\Phi:\A \rightarrow M_3(\R)$ defined by $\Phi(f) = \rho_w(f) \oplus \Pi_w(f)$, for some non-symmetric point $w \in \Omega$. By Lemma \cite[Lemma 4.2]{ANSS}, one can 
choose $g \in \A$ such that $g(w) = i$. 
Observe that $\Phi(g)^2 = 0 \, \oplus \, -I_2$. 
Suppose that $\Phi(g)^2 \in \operatorname{Range}(\Phi)$. Say, $\Phi(h) = \Phi(g)^2$. It must be 
the case that $Re \, h(w) = 0$ and $h(w) = -1$ simultaneously, which is impossible. Note also that 
$\Phi$ is not a homomorphism, since $\rho_w$ is not multiplicative due to the fact that 
$w \neq \bar w$. Thus in general, a $C^*$-extreme UCP map $\Phi:\A \rightarrow M_n(\R)$
 need not be multiplicative and its range need not be an algebra. Note that this is in stark contrast to the nature of C$^*$-extreme UCP maps on commutative unital complex C$^*$-algebras, which necessarily are $*$-homomorphisms. 
 \vspace{0.2 cm}

(ii) Consider  the linear extreme UCP map $\Phi_S:\C \rightarrow M_3(\R)$ defined 
by $\Phi_S(\lambda) = 
Re(\lambda)I_3 + Im(\lambda)S$, where $S = \begin{pmatrix}  
0 & \alpha & \beta \\ -\alpha & 0 & \gamma \\ - \beta & -\gamma & 0
\end{pmatrix}$ with $\|S\| = \alpha^2 + \beta^2 + \gamma^2 = 1$.  
It was observed in \cite[example 4.16]{ANSS} that $\Phi_S$ is not 
$C^*$-extreme except when $(\alpha, \beta, \gamma) = (1, 0, 0)$. 
But it actually is, by Theorem \ref{thm:range}, since 
$\operatorname{Range}(\Phi_S)$ is commutative. \\ 
\end{remark}

\section{Domain Classification}
\label{sec:main2} 
This section contains an application of our main result. i.e. Theorem 
\ref{thm:range}. It is divided into three subsections containing several 
propositions that lead us finally to a proof of Theorem \ref{thm:dom}.  
\subsection{The model algebras} 
\label{subsec:models}
Here, we establish that on the fundamental real unital $C^*$-algebras $\R$, $\R^2$ and 
$\C$, the notions of linear extremality and $C^*$-extremality coincide for UCP maps. 

\begin{proposition}
\label{prop:models} 
Let $n \ge 1$. If $\mathcal A$ is either $\R$ or $\R^2$ or $\C$, then 
$\Phi \in \operatorname{UCP}(\mathcal A, M_n(\R))$ is linear extreme 
if and only if it is $C^*$-extreme. 
\end{proposition}

\begin{proof}
In view of Remark \ref{rem:oneway}, 
it suffices to show that linear extremality implies $C^*$-extremality. \\

If $\A = \R$, then there exists only one unital linear map $\Phi:\A \rightarrow M_n(\R)$, which 
is given by $\Phi(a) = a I_n$. Since it is a $*$-homomorphism, it is both $C^*$-extreme and linear 
extreme. \\

If $\A = \R^2$, then any UCP map $\Phi:\A \rightarrow M_n(\R)$ is of the form 
$\Phi(a,b) = a E + b (I_n - E)$, where $E = \Phi(1,0)$ and $I_n - E = \Phi(0,1)$. 
Since $(1,0) \text{ and } (0,1)$ are positive elements in $\A$, it follows that 
$0 \le E, I_n - E \le I_n$. By (the real version of) \cite[Corollary 1.4.4]{A}, the UCP map $\Phi$
is linear extreme if and only if $E$ and $I_n - E$ are projections. In particular, the 
range of the linear extreme UCP map $\Phi$ is commutative.  
It follows from Theorem \ref{thm:range} that $\Phi$ is also $C^*$-extreme. \\

If $\A = \C$, by \cite[Proposition 4.13]{ANSS}, any UCP map 
$\Phi:\A \rightarrow M_n(\R)$ is of the form 
$\Phi(a + ib) = a I_n + b S$, where $S = \Phi(i) \in M_n(\R)$ is a skew-symmetric matrix, i.e., 
$S^* = -S$ satisfying $\|S\|\le 1$. Note that the UCP map $\Phi$ has commutative range. It follows from Theorem \ref{thm:range} that every linear extreme UCP map is also $C^*$-extreme. 
\end{proof}

\begin{remark}
\label{rem:ext}
Combining the conclusion of Proposition \ref{prop:models} for $\A = \C$, with Proposition \cite[Proposition 4.13]{ANSS}, one gets that for the set 
\[
 \mathcal {S}_n = \{ T \in M_n(\R)\,:\, T^* = -T \text{ and } \|T\|\le 1 \} 
 \]
 of $n \times n$ skew-symmetric contractions, the (linear) extreme points and the $C^*$-extreme points coincide. 
\end{remark}

\begin{remark}
\label{rem:non-models}
Recall $\A = C(\Omega, -)$.
If $\A \not \cong \{\R, \R^2, \C\}$, equivalently, if $\dim(\A) \ge 3$, then one of the following must hold:
\begin{enumerate}
\item[(i)] $\Omega$ contains at least three symmetric points; or
\item[(ii)] $\Omega$ contains at least one symmetric point and one non-symmetric point; or
\item[(iii)] $\Omega$ contains at least two distinct non-symmetric points that are not conjugates of each other. 
\end{enumerate}
By Lemma \cite[Lemma 4.2]{ANSS}, according as (i) or (ii) or (iii) holds, there exists a surjective $*$-homomorphism from $\A$ onto $\R^3$ or $\R \oplus \C$, or $\C \oplus \C$, respectively.  
\end{remark}

\subsection{The violator algebras}
\label{subsec:violators}
In view of Remark \ref{rem:non-models}, in this subsection 
we consider the real unital $C^*$-algebras $\R^3$, 
$\R \oplus \C$ and $\C \oplus \C$ and explore the coincidence of 
the notions of linear and $C^*$-extremality for UCP maps on them.  

\begin{lemma}
\label{lem:R3}
For every $n \ge 2$, there exists a linear extreme UCP map 
$\Psi:\R^3 \rightarrow M_n(\R)$ that is not $C^*$-extreme. 
\end{lemma} 

\begin{proof}
First assume $n = 2$. The result \cite[Theorem 1.4.10]{A} 
gives the explicit structure
of linear extreme UCP maps on unital commutative C$^*$-algebras. 
The following example that is well-known in the context 
of quantum measurements (see for instance \cite[Section 2.1.4]{C}) 
is one such UCP map, but adapted to our real setting. 
Let 
\[
u_1=\begin{pmatrix}  1 \\  0 \end{pmatrix},
\qquad
u_2=\begin{pmatrix} -\frac12 \\ \frac {\sqrt{3}}{2} \end{pmatrix}, 
\qquad
u_3=\begin{pmatrix} - \frac 12 \\ -\frac {\sqrt{3}}{2} \end{pmatrix}, 
\qquad
E_i=\frac23u_iu_i^*.
\]
Define the linear map $\Gamma:\R^3 \rightarrow M_2(\R)$ by 
\[
\Gamma(x_1,x_2,x_3)=x_1E_1+x_2E_2+x_3E_3. 
\]
Note that $\Gamma(x_1, x_2, x_3) 
= V^* \pi(x_1, x_2, x_3) V$, where  $\pi:\R^3 \rightarrow M_3(\R)$ is the 
$*$-homomorphism and $V:\R^2 \rightarrow \R^3$ is the isometry given by 
\[
\pi(x_1,x_2,x_3)= \begin{pmatrix}  x_1 & 0 & 0 \\ 0 & x_2 & 0 \\ 0 & 0 & x_3 \end{pmatrix} \text{ and } 
V=\sqrt{2/3} \begin{pmatrix} u_1^*\\u_2^*\\u_3^* \end{pmatrix}.
\]
The linear extremality of $\Gamma$  follows from \cite[Theorem 1.4.10]{A} and 
Theorem \ref{thm:arveson}. But 
we verify it directly here.  First, observe that $(\R^3, \pi, V)$ is a minimal Stinespring representation for 
$\Gamma$. This is immediate from the fact that
 $h = u_1$, $\pi(x_1, x_2, x_3)Vh = \sqrt{\frac 23} (x_1, -\frac12 x_2, -\frac12 x_3)$. 
Next, note that $\pi(\R^3)' = \pi(\R^3)'_{sa} = \pi(\R^3)$. In particular,  
every element in the commutant is self-adjoint. 
Let $T = \begin{pmatrix}  t_1 & 0 & 0 \\ 0 & t_2 & 0 \\ 0 & 0 & t_3 \end{pmatrix} \in \pi(\R^3)'_{sa}$ be 
such that $V^*TV= \frac 2 3 \sum_{i=1}^3 t_iE_i = 0.$
Since $\{E_1,E_2,E_3\} \subset M_2(\R)$ is linearly independent, it follows that 
$t_1 = t_2 = t_3 = 0$. Thus $T = 0$. By Theorem \ref{thm:arveson}, the UCP map 
$\Gamma$ is linear extreme.  Note that $\Gamma(1,0,0) = E_1$ and 
$\Gamma(0,1,0) = E_2$, do not commute. \\

Suppose $n \ge 3$, define the UCP map $\Psi: \R^3 \rightarrow M_n(\R)$ by  
$
\Psi = \Gamma \oplus \rho^{(n-2)}, 
$
where $\rho:\R^3 \rightarrow \R$ is defined by $\rho(x_1, x_2, x_3) = x_1$. Note 
that $\rho^{(n-2)}$ is a $*$-homomorphism. By Lemma \ref{lem:dirsum}, 
it follows that $\Psi$ is a linear extreme UCP map.  Note that the range of 
$\Psi$ is not commutative due to the noncommutativity of the
 range of $\Gamma$. By Theorem \ref{thm:range}, the UCP map 
$\Psi$ is not $C^*$-extreme. 
\end{proof}

\begin{lemma}
\label{lem:RC}
For every $n \ge 2$, there exists a linear extreme UCP map  
$\Psi:\R \oplus \C \rightarrow M_n(\R)$ that is not $C^*$-extreme.    
\end{lemma}

\begin{proof}
Assume $n = 2$. 
Let $\rho:\R \oplus \C \rightarrow \R$ and $\zeta:\R \oplus \C \rightarrow M_2(\R)$ be 
defined by 
\[
\rho(a,z)=a,
\quad
\zeta(a,z)=
\begin{pmatrix}
Re \, z & Im \, z\\
-Im \, z & Re \, z
\end{pmatrix}.
\]
Fix $0<t<1$, and define 
\begin{equation*}
\begin{split}
&A=\begin{pmatrix}\sqrt t &0\end{pmatrix}:
\R^2 \to\R,  
\,
B =\begin{pmatrix}
\sqrt{1-t}&0 \\
0&1
\end{pmatrix}:
\R^2\to\R^2 \\
&\text{and } V = \begin{pmatrix} A \\ B \end{pmatrix}:\R^2 \rightarrow \R^3.
\end{split}
\end{equation*}
Note that $V$ is an isometry. The definition of this isometry $V$ is inspired by 
\cite[Lemma 4.3]{BH}. Let $\pi:\R \oplus \C \rightarrow M_3(\R)$ be the 
unital $*$-homomorphism 
defined by $\pi = \rho \oplus \zeta$ and $\Gamma: \R \oplus \C \rightarrow M_2(\R)$ be the 
UCP map defined by 
$
\Gamma(a,z) =V^*\pi(a,z)V.
$
Since the maps $A$ and $B$ are surjective, it follows that 
 the Stinespring representation $(\R^3, \pi, V)$ of the UCP map $\Gamma$
 is minimal. Let $T \in \pi(\R\oplus\C)'_{sa}=\R\oplus \R I_2$. 
Say $T = r \oplus s I_2$ for some $r, s \in \R$. If 
$V^*TV = rA^*A + s B^*B = 0$, then evaluating at standard basis 
$e_1, e_2$ of $\R^2$ yields, $r = s = 0$. Thus $T = 0$. By 
Theorem \ref{thm:arveson}, it follows that $\Gamma$ is a linear extreme 
UCP map. Note that $\Gamma(1,0) = t \begin{pmatrix}  1 & 0 \\ 0 & 0 \end{pmatrix}$ 
and $\Gamma(0, i) = \sqrt{1 - t} \begin{pmatrix}  0 & 1  \\ -1 & 0 \end{pmatrix}$ 
do not commute. \\

Suppose $n \ge 3$, define the UCP map $\Psi:\R \oplus \C \rightarrow M_n(\R)$ by 
$\Psi =\Gamma \oplus \rho^{(n-2)}$. By Lemma \ref{lem:dirsum}, 
$\Psi$ is linear extreme. Note that the range of 
$\Psi$ is not commutative due to the noncommutativity of the
 range of $\Gamma$. By Theorem \ref{thm:range}, 
$\Psi$ is not $C^*$-extreme. 
\end{proof}

\begin{lemma}
\label{lem:CC}
For every $n \ge 2$, there exists a linear extreme UCP map 
$\Psi:\C \oplus \C \rightarrow M_n(\R)$ that is not $C^*$-extreme. 
\end{lemma}

\begin{proof}
First assume $n = 2$. Fix $t \in (0,1)\setminus \{1/2\}$, 
 and define $A,B: \R^2 \rightarrow \R^2$ and 
$V:\R^2 \rightarrow \R^4$ by 
\begin{equation*}
A=\begin{pmatrix} \sqrt t &0\\0& \sqrt {1 - t} \end{pmatrix}, 
\quad
B=\begin{pmatrix} \sqrt{1- t} &0\\0& \sqrt{t} \end{pmatrix}
\text{ and } 
V = \begin{pmatrix}  A \\ B \end{pmatrix}. 
\end{equation*}
Note that $V$ is an isometry. The definition of this isometry $V$ is inspired by 
\cite[Lemma 4.3]{BH}.\\

Let $\pi: \C \oplus \C \rightarrow M_4(\R)$ denote the unital $*$-homomorphism defined by 
\[
\pi(z,w) = \begin{pmatrix} Re\, z & Im\, z \\ - Im\, z & Re\, z \end{pmatrix} 
\oplus \begin{pmatrix} Re\, w & Im\, w \\ - Im\, w & Re\, w \end{pmatrix},
\]
and 
$\Gamma: \C \oplus \C \rightarrow M_2(\R)$ denote the UCP map defined by 
$\Gamma(z,w) =V^*\pi(z,w)V.$  Minimality of the Stinespring representation 
$(\R^4, \pi, V)$ for the UCP map $\Gamma$ follows from the surjectivity of 
$A$ and $B$. Note that $\pi(\C \oplus \C)'_{sa} = \R I_2 \oplus \R I_2$. Let 
$T \in \pi(\C \oplus \C)'_{sa}$. Say $T = 
rI_2 \oplus sI_2$ for some $r, s \in \R$. If $V^*TV = rA^*A+ sB^*B = 0$, then it follows 
from a simple computation that $T = 0$. 
Moreover, 
\[
\Gamma(1,0)= \begin{pmatrix} t & 0 \\ 0 & 1 - t \end{pmatrix} 
\text{  and  } \Gamma(i,0) = \sqrt{t (1 - t)} \begin{pmatrix} 0 & 1 \\ - 1 & 0 
\end{pmatrix}
\]
do not commute. \\

Assume $n = 3$. Fix $t \in (0,1) \setminus \{1/2\}$ and define the linear maps 
 $A, B, C: \R^3 \rightarrow \R^2$ and $V:\R^3 \to \R^6$ by 
 \[
 A=\begin{pmatrix} \sqrt t & 0& 0\\ 0 & \sqrt {1 - t} & 0\end{pmatrix}, 
 \quad
 B=\begin{pmatrix} \sqrt{1- t} &0 & 0 \\0& \sqrt{t} & 0 \end{pmatrix},
 \quad
 C = \begin{pmatrix} 0 & 0 & 1 \\0 & 0 & 0\end{pmatrix}
  \]
 $
\text{ and }
\quad  
 V = \begin{pmatrix}  A \\ B \\ C\end{pmatrix}. 
 $ 
 Note that $V$ is an isometry. Let $\pi_1:\C \oplus \C \rightarrow M_2(\R)$ be the unital $*$-homomorphism 
defined by 
\[
\pi_1(z,w) = \begin{pmatrix} Re\, z & Im\, z \\ - Im\, z & Re\, z \end{pmatrix}. 
\]
Define the UCP map
 $\Delta:\C \oplus \C  \rightarrow M_3(\R)$ by 
$\Delta(z,w) =V^*(\pi \oplus \pi_1)(z,w)V$. For the standard basis vector $e_3 \in \R^3$, 
observe that 
\[
\{(\pi \oplus \pi_1)(z,w)Ve_3\,:\,(z,w) \in \C \oplus \C\} = \{(Re \, z, - Im \, z)\,:\, z \in \C\} = \R^2.
\]
This, along with the fact that $A$ and $B$ are surjective implies that 
$(\R^6, \pi \oplus \pi_1, V)$ is a minimal Stinespring dilation for the UCP map $\Delta$.
\\

 To prove linear extremality of $\Delta$, let $T \in (\pi \oplus \pi_1)(\C\oplus \C)'_{sa}$. Note that 
 $T = \begin{pmatrix} 
 \alpha I_2 & 0 & x \\ 0 & \beta I_2 & 0 \\x^* & 0 & \gamma I_2 \end{pmatrix}$, where $\alpha, \beta, \gamma \in \R$ 
and $x = \begin{pmatrix} a & b \\ c & d \end{pmatrix} \in M_2(\R)$. 
Compute 
 \[
  V^*TV = \begin{pmatrix}  \alpha t + \beta (1 - t) & 0 & a \sqrt{t} \\
        0 & \beta t + \alpha (1-t) & c \sqrt{1-t} \\
        a \sqrt{t} & c \sqrt{1-t} & \gamma         
        \end{pmatrix}. 
  \]
If $V^*TV = 0$, then $c = a = \gamma  = 0$, since $t > 0$. Also since $2t - 1 \neq 0$, 
$\alpha t + \beta (1 - t) = \beta t + \alpha (1-t) = 0$ forces $\alpha = \beta = 0$. 
Finally, since $T \in (\pi \oplus \pi_1)(\C\oplus \C)'_{sa}$, it must be the case that the 
matrix $x$ commutes with $\begin{pmatrix} 0 & 1 \\ -1 & 0 \end{pmatrix}$. 
This forces $T = 0$. It follows from Theorem \ref{thm:arveson} that $\Delta$ is a linear 
extreme UCP map. Moreover, 
\[
\Delta(1,0)= \begin{pmatrix} t & 0 & 0\\ 0 & 1 - t & 0 \\0 & 0 & 1\end{pmatrix} \text{  and  }
\Delta(i,0) = \sqrt{t (1 - t)} \begin{pmatrix} 0 & 1 & 0 \\ - 1 & 0 & 0 \\ 0 & 0 & 0
\end{pmatrix}
\]
do not commute. \\


Let $n\ge 4$. 
If $n=2k$ is even, define
\[
\Psi :=\Gamma \oplus \pi_1^{(k-1)}.
\]

If $n=2k+1$ is odd, define 
\[
\Psi:=\Delta \oplus \pi_1^{(k-1)}. 
\]
In both cases, $\Psi:\C \oplus \C \rightarrow M_n(\R)$ is a linear 
extreme UCP map by Lemma \ref{lem:dirsum}. Since $\Gamma$ and 
$\Delta$ have non-commutative range, it follows from Theorem 
\ref{thm:range} that $\Psi$ is not a $C^*$-extreme UCP map.  
\end{proof}

\subsection{Proof of Theorem \ref{thm:dom}}

We are now in a position to prove the second main result of 
the paper, i.e., Theorem \ref{thm:dom}.

\begin{proof}
Recall the unital commutative $C^*$-algebra $\A = C(\Omega, -).$ 
The equivalence $(ii) \iff (iii)$ is 
immediate. \\

To prove $(i) \implies (ii)$, suppose that $\mathcal A$ is non-isomorphic to $\R$ 
or $\R^2 $ or $\C$. 
By Remark \ref{rem:non-models}, there exists a surjective $*$-homomorphism $\rho: \mathcal A \rightarrow 
\mathcal B$, where $\mathcal B \in \{\R^3, \R \oplus \C, \C \oplus \C\}$. 
Let $\Psi:\mathcal B \rightarrow M_n(\R)$ be the linear extreme but not $C^*$-extreme 
UCP map in Lemmas \ref{lem:R3}, \ref{lem:RC} and \ref{lem:CC}, according to 
the choice of $\mathcal B$. Note in particular that $\Psi$ does not have 
commutative range. By Lemma \ref{lem:comp} it follows that the UCP map 
$\Phi = \Psi \circ \rho: \mathcal A \rightarrow M_n(\R)$ is linear extreme. But since 
$\rho$ is surjective, the range of $\Phi$ equals the range of $\Psi$, which is not 
commutative. Hence $\Phi$ is not $C^*$-extreme, by Theorem \ref{thm:range}.
 This is a contradiction. \\

To prove $(ii) \implies (i)$, Let $\mathcal C \in \{\R, \R^2,\C\}$ and 
let $\rho:\mathcal C \rightarrow \mathcal A$ be a $*$-isomorphism. If $\Phi:\mathcal A 
\rightarrow M_n(\R)$ is a linear extreme UCP map, then it follows from 
Lemma \ref{lem:comp} that $\Phi \circ \rho:\mathcal C \rightarrow M_n(\R)$ is 
also a linear extreme UCP map. 
Suppose that $\Phi$ is expressed as a proper $C^*$-convex combination
of the UCP maps $\Phi_1,\dots,\Phi_k:\mathcal A \rightarrow M_n(\R)$, 
say $\Phi(a) = \sum_{i=1}^k T_i^* \Phi_i(a) T_i$, for each $a \in \A$. 
Since $\rho$ is bijective, 
$
(\Phi \circ \rho)(x) = \sum_{i=1}^k T_i^* ((\Phi_i \circ \rho)(x)) T_i,
$
for each $x \in \mathcal C$. By Proposition \ref{prop:models}, the linear extreme 
UCP map $\Phi \circ \rho:\mathcal C \rightarrow M_n(\R)$ is also $C^*$-extreme. 
Thus there exists unitaries $U_i \in M_n(\R)$ such that 
$U_i^* \Phi_i(\rho(x)) U_i = \Phi(\rho(x))$ for each $x \in \mathcal C$. i.e., 
$\Phi(a) = U_i^* \Phi_i(a)U_i$ for some unitaries $U_i \in M_n(\R)$ and 
each $a \in \A$. Hence $\Phi$ is a $C^*$-extreme UCP map, and 
the proof is complete. 
\end{proof}
\noindent
{\bf Acknowledgements:} The authors thank Professor Scott McCullough (UF) for   
his valuable feedback and suggestions.  \\ 

\noindent
{\bf Statement on AI use:}  The authors used ChatGPT 
for preliminary literature searches, identifying relevant references, checking certain computations 
and proof-reading. During proof-reading, the AI system  flagged two minor errors in a computation
 concerning a commutant in Lemma \ref{lem:CC}, which the authors subsequently verified and 
corrected. The mathematical strategy, ideas, proofs and results presented in this article were 
developed by the authors, who take full responsibility for their mathematical correctness. 



\begin{thebibliography}{99}
\bibitem[A]{A} W.~B. Arveson, \emph{Subalgebras of $C^*$-algebras}, 
Acta Math. \textbf{123} (1969), 141--224.

\bibitem[ANSS]{ANSS}
Anand O. R., K. Sumesh and A. Sutradhar,
\emph{$C^*$-extreme points of unital completely positive maps on real
$C^*$-algebras}, arXiv:2502.15362.

\bibitem[ANSS1]{ANSS1} Anand O. R., K. Sumesh, and A. Sutradhar,
\emph{$C^*$-extreme points of unital completely positive maps on real
$C^*$-algebras}, Studia Math. \textbf{288} (2026), 201--232. 

\bibitem[BH]{BH} S.~Balasubramanian and N.~Hotwani, 
\emph{$C^*$-extreme entanglement breaking maps on operator systems}, 
Linear Algebra Appl. \textbf{685} (2024), 182--213.

\bibitem[BBK]{BBK} T.~Banerjee, B.~V. Rajarama Bhat, and M.~Kumar,
\emph{$C^*$-extreme points of positive operator valued measures and unital completely positive maps},
 Comm. Math. Phys. \textbf{388} (2021), no.~3, 1235--1280. 

\bibitem[BK]{BK} B.~V. Rajarama Bhat and M.~Kumar, 
\emph{$C^*$-extreme maps and nests}, 
J. Funct. Anal. \textbf{282} (2022), no.~8, Paper No.~109397.

\bibitem[BDMS]{BDMS} B.~V. Rajarama Bhat, R.~Devendra, N.~Mallick, and K.~Sumesh, 
\emph{$C^*$-extreme points of entanglement breaking maps}, 
Rev. Math. Phys. \textbf{35} (2023), no.~3, Paper No.~2350005.

\bibitem[C]{C} R. Cleve, \emph{Notes on Quantum Information Theory}, Nov 2023,  
\url{https://cleve.iqc.uwaterloo.ca/resources/QIC-710-F23/Qic710QuantumInformationTheory2023V2.pdf}

\bibitem[F]{F} D.~R. Farenick, \emph{Extremal matrix states on operator systems},
 J. London Math. Soc. (2) \textbf{61} (2000), no.~3, 885--892. 

\bibitem[FM1]{FM1}
D. R. Farenick and P. B. Morenz,
\emph{$C^*$-extreme points in the generalised state spaces of a
$C^*$-algebra},
Trans. Amer. Math. Soc. 349 (1997), 1725--1748.

\bibitem[FM2]{FM2} D.~R. Farenick and P.~B. Morenz,
 \emph{$C^*$-extreme points of some compact $C^*$-convex sets}, 
Proc. Amer. Math. Soc. \textbf{118} (1993), no.~3, 765--775.

\bibitem[FZ]{FZ}
D. R. Farenick and H. Zhou,
\emph{The structure of $C^*$-extreme points in spaces of completely
positive linear maps on $C^*$-algebras},
Proc. Amer. Math. Soc. 126 (1998), 1467--1477.

\bibitem[G1]{G1} M.~C. Gregg,
 \emph{$C^*$-Extreme Points of the Generalized State Space of a Commutative $C^*$-Algebra}, 
Ph.D. thesis, The University of Nebraska--Lincoln, ProQuest LLC, Ann Arbor, MI, 2008. 

\bibitem[G2]{G2} M.~C. Gregg, 
\emph{On $C^*$-extreme maps and $*$-homomorphisms of a commutative $C^*$-algebra},
Integral Equations Operator Theory \textbf{63} (2009), no.~3, 337--349.

\bibitem[HJ]{HJ} R. A. Horn and C. R. Johnson, {\it Matrix Analysis}, 2nd edition, Cambridge University Press,
 New York, 2012. 

\bibitem[HMP]{HMP} A.~Hopenwasser, R.~L. Moore, and V.~I. Paulsen,
 \emph{$C^*$-extreme points}, Trans. Amer. Math. Soc. \textbf{266} (1981), no.~1, 291--307.

\bibitem[L]{L} B.~Li, \emph{Real Operator Algebras}, World Scientific, River Edge, NJ, 2003. 


\bibitem[LP]{LP} R. I. Loebl and V. I. Paulsen,
\emph{Some remarks on $C^*$-convexity},
Linear Algebra Appl. 35 (1981), 63--78.

\bibitem[LT]{LT} C. K. Li and N. K. Tsing, \emph{Duality between some linear preserver
problems. III. c-spectral norms on (skew)-symmetric matrices and matrices with fixed singular values,} Linear Algebra Appl. 143:67-97 (1991).


\bibitem[P]{P} Paulsen, V. I., \emph{Completely bounded maps and operator algebras,} vol. 78, 
Cambridge Studies in Advanced Mathematics, Cambridge University Press, 2002. ISBN 0-521-81669-6.

\bibitem[S]{S} W.~F. Stinespring,
\emph{Positive functions on $C^*$-algebras}, 
Proc. Amer. Math. Soc. \textbf{6} (1955), 211--216.

\bibitem[Z]{Z} H.~Zhou, 
\emph{$C^*$-extreme points in spaces of completely positive maps},
 Ph.D. thesis, University of Regina, 1998.

\end{thebibliography}
\end{document}